\documentclass[11pt]{article}
\usepackage[english]{babel}
\usepackage[letterpaper,top=2cm,bottom=2cm,left=3cm,right=3cm,marginparwidth=1.75cm]{geometry}

\usepackage{amsmath}
\usepackage{graphicx}
\usepackage{amssymb}
\usepackage{algorithm2e}
\usepackage[colorlinks=true, allcolors=blue]{hyperref}
\usepackage{color}
\usepackage{amsthm}

\newtheorem*{theorem}{Theorem}

\title{A generalization of Boole's formula derived from a system of linear equations}
\author{Haoran Zhu \thanks{College of Sciences, Northeastern University, Shenyang, Liaoning 110004, China ({\tt whrzhu@outlook.com}).}
}

\begin{document}
	\maketitle
	
	\begin{abstract}
		We analyze a system of linear algebraic equations whose solutions lead to a proof of a generalization of Boole's formula. In particular, our approach provides an elementary and short alternative to Katsuura's proof of this generalization.
	\end{abstract}
	
	
	\section{Introduction}  
	Due to their numerous relations, binomial coefficients play an important role in various mathematical fields, including 
	enumerative combinatorics, statistics and number theory. 
	In Boole's classical book “Calculus of Finite Differences” \cite{Bo}, the following beautiful formula is given, which holds for $1\le m\le n \in \mathbb N$:
	
	\begin{align}\label{boole}
		\sum_{k=0}^{n}(-1)^{n-k}\binom{n}{k}k^m=
		\begin{cases} 
			n! & \text{if }  m = n, \\
			0  & \text{if }  m < n.
		\end{cases}
	\end{align}
	
	It is known that the formula is related to Stirling's partition numbers $S(m, n)$ (see, e.g.,\cite{HH}), which is given by the following equation
	\[
	\sum_{k=0}^{n}(-1)^{n-k}\binom{n}{k}k^m=n!\cdot S(m,n) \quad (m,n\in \mathbb{N}).
	\]
	The relations between the formulas have implications for the study of the degrees of normed null-polynomials and the derivation of inequalities associated with the Smarandache function (see \cite{SHW1,SHW2}).\par
	
	The enduring interest in Boole's formula has led to a variety of proof techniques being developed over the years. Gould \cite{Gou1} discussed its properties and called it {\em Euler's formula}. In 2005, Anglani and Barile \cite{Ab} introduced two proofs via methods from real analysis and combinatorics. Subsequently, Phoata \cite{P} and Katsuura \cite{K} provided new proofs and gave a generalization of Boole's formula. In this note, we give a short and elementary proof of this formula which is based on a system of linear equations. More recently, Alzey and Chapman \cite{Ac} have presented a novel proof, while Batır and Atpınar \cite{Ba,BA} have independently developed two entirely new approaches to validating Boole's formula.\par
	
	\section{Solutions of linear algebraic equations}
	Let us consider the following case
	\begin{equation}\label{eq}
		\begin{bmatrix}
			1 & 1 & \cdots & 1 \\
			a & a+b & \cdots & a+nb \\
			a^2 & (a+b)^2 & \cdots & (a+nb)^2 \\
			\vdots & \vdots & \ddots & \vdots \\
			a^n & (a+b)^n & \cdots & (a+nb)^n
		\end{bmatrix}
		\begin{bmatrix}
			x_0 \\
			x_1 \\
			x_2 \\
			\vdots \\
			x_n
		\end{bmatrix} = \begin{bmatrix}
			0 \\
			0 \\
			0 \\
			\vdots \\
			b^n \cdot n!
		\end{bmatrix}.
	\end{equation}
	where $a,b$ can be any real numbers and the coefficient matrix $V$ is a Vandermonde matrix.\par 
	To solve this system, we first calculate the determinant of $V$. Since $V$ is a Vandermonde matrix, its determinant can be computed as follows:
	\begin{equation*}\label{det_V}
		\det(V) = \prod_{0 \leq j < i \leq n} ((a+ib) - (a+jb)) = n!\cdot (n-1)!\cdots 1!\cdot b^{\frac{n(n+1)}{2}}.
	\end{equation*}\par
	We then proceed to define matrix $V_k$ as follows:
	\[
	V_k = 
	\begin{bmatrix}
		1 & 1 & \cdots & 1 & 0 & 1 & \cdots & 1 \\
		a & a+b & \cdots & a+(k-1)b & 0 & a+(k+1)b & \cdots & a+nb \\
		\vdots & \vdots & \ddots & \vdots & \vdots & \vdots & \ddots & \vdots \\
		a^n & (a+b)^n & \cdots & (a+(k-1)b)^n & b^n \cdot n! & (a+(k+1)b)^n & \cdots & (a+nb)^n
	\end{bmatrix}
	\]
	and we denote the matrix obtained by removing the $(n+1)$-th row and $(k+1)$-th column from $V_k$ as $V_{k1}$. The determinant of $V_k$ is computed by first applying Laplace's expansion along the $(k+1)$-th column, which yields the determinant of the submatrix $V_{k1}$. This submatrix is also a Vandermonde matrix, and its determinant can be computed using the way previously showed.
	\begin{align*}
		\det(V_k) &= (-1)^{n-k}\cdot b^n\cdot n!\cdot \det(V_{k1})\\ 
		&= (-1)^{n-k}\cdot b^{\frac{n(n+1)}{2}}\cdot n!\cdot \frac{n!\cdots (k+1)!\cdot (k-1)!\cdots 1!}{(n-k)!}.
	\end{align*}\par
	
	It is sufficient for us to use Cramer's rule to obtain an explicit expression for $x_k$
	\begin{align*}
		x_k=\frac{\det(V_k)}{\det(V)}=\frac{(-1)^{n-k}\cdot n!}{k!\cdot (n-k)!}=(-1)^{n-k}\binom{n}{k}.
	\end{align*}
	
	Reading the equations in (\ref{eq}) one by one we obtain the following result.
	\begin{theorem}(Generalization of Boole's formula)
		For any real numbers $a$ and $b$, and for $1\le m\le n \in \mathbb N$, we have
		\begin{equation}\label{geq}
			\sum_{k=0}^{n}(-1)^{k}\binom{n}{k}(a+bk)^m=
			\begin{cases} 
				(-1)^{n}\cdot b^n\cdot n! & \text{if }  m = n, \\
				0  & \text{if }  m < n.
			\end{cases}
		\end{equation}
	\end{theorem}

	For the special case $a = 0$ and $b = 1$, the result in the Theorem implies Boole's formula (\ref{boole}).
	
	Perhaps similar approaches from linear algebra can also be used to generalize other combinatorial identities.
	
	\bigskip
	\noindent{\bf Acknowledgements
	}The author thanks the referees for their constructive feedback, which enhanced the paper's organization.
	
	\bibliographystyle{plain}

\end{document}